\documentclass[12pt]{amsart}
\usepackage{amssymb}
\usepackage{graphicx}

\newtheorem{Def}{Definition}
\newtheorem{Lem}{Lemma}

\newtheorem{Thm}{Theorem}
\newtheorem{Cor}{Corollary}

\newenvironment{Pf}{ Proof.}{\(\square\)}
\newtheorem{Exm}{Example}

\title[On generalized Berwald surfaces with locally symmetric fourth root metrics]{On generalized Berwald surfaces with locally symmetric fourth root metrics}
\author{Cs. Vincze, T. Khoshdani and M. Ol\'{a}h}
\address{Inst. of Math., Univ. of Debrecen \\
H-4010 Debrecen, P.O.Box 12 \\
Hungary}
\email{csvincze@science.unideb.hu}
\email{khoshdani@yahoo.com, olma4000@gmail.com}
\begin{document}
\begin{abstract} Let $m=2l$ be a positive natural number, $l=1, 2, \ldots. $ A Finslerian metric $F$ is called an $m$-th root metric if its $m$-th power $F^m$ is of class $C^{m}$ on the tangent manifold $TM$. Using some homogenity properties, the local expression of an $m$-th root metric is a polynomial of degree $m$ in the variables $y^1$, $\ldots$, $y^n$, where $\dim M=n$. $F$ is locally symmetric if each point has a coordinate neighbourhood such that $F^m$ is a symmetric polynomial of degree $m$ in the variables $y^1$, $\ldots$, $y^n$ of the induced coordinate system on the tangent manifold. Using the fundamental theorem of symmetric polynomials, the reduction of the number of the coefficients depending on the position makes the computational processes more effective and simple. In the paper we present some general observations about locally symmetric  $m$-th root metrics. Especially, we are interested in generalized Berwald surfaces with locally symmetric fourth root metrics. The main result (Theorem 1) is their intrinsic characterization in terms of the basic notions of linear algebra. We present a one-parameter family of examples as well. The last section contains some computations in 3D. They are supported by the MAPLE mathematics softwer (LinearAlgebra).
\end{abstract}
\maketitle
\footnotetext[1]{Cs. Vincze is supported by the EFOP-3.6.1-16-2016-00022 project. The project is co-financed by the European Union and the European Social Fund.}
\footnotetext[2]{T. Khoshdani is supported by the Department of Mathematics of University of Mohaghegh Ardabili, Ardabil, Iran.}
\section*{Introduction}
Let $M$ be a differentiable manifold with local coordinates $u^1, \ldots, u^n.$ The induced coordinate system of the tangent manifold $TM$ consists of the functions $x^1, \ldots, x^n$ and $y^1, \ldots, y^n$. For any $v\in T_pM$, $x^i(v)=u^i\circ \pi(v)$ and $y^i(v)=v(u^i)$, where  $\pi \colon TM\to M$ is the canonical projection, $i=1, \ldots, n$. Introducing the so-called Liouville vector field $\displaystyle{C:=y^i\partial/\partial y^i}$, it can be easily seen that the integral curves of $C$ are of the form $e^tv$ for any nonzero $v\in TM$. If the function $f\colon TM\to \mathbb{R}$ is differentiable on the complement of the zero section such that $f(tv)=t^kf(v)$ for any positive $t\in \mathbb{R}$, then we have that
$$C(v)f=\frac{d}{dt}\left(f(e^tv)\right)_{t=0}=kf(v)$$
and vice versa. This means that the function $f$ is positively homogeneous of degree $k$ if and only if $Cf=kf$. It is the so-called Euler's theorem on homogeneous functions. 

A Finsler metric \cite{BSC} is a continuous function $F\colon TM\to \mathbb{R}$ satisfying the following conditions: $\displaystyle{F}$ is smooth on the complement of the zero section (regularity), $\displaystyle{F(tv)=tF(v)}$ for all $\displaystyle{t> 0}$ (positive homogenity) and the Hessian $\displaystyle{g_{ij}=\partial^2 E /\partial y^i \partial y^j}$, where $E=\frac{1}{2}F^2$, is positive definite at all nonzero elements $\displaystyle{v\in TM}$ (strong convexity). 

A linear connection $\nabla$ on the base manifold $M$ is called \emph{compatible} to the Finslerian metric if the parallel transports with respect to $\nabla$ preserve the Finslerian length of tangent vectors. Finsler manifolds admitting compatible linear connections are called generalized Berwald manifolds. It can be easily seen \cite{VKMO} that a linear connection  $\nabla$ on the base manifold $M$ is compatible to the Finslerian metric function if and only if the induced horizontal distribution is conservative, i.e. the derivatives of the fundamental function $F$ vanish along the horizontal directions with respect to $\nabla$:
\begin{equation}
\label{cond1}
\frac{\partial F}{\partial x^i}-y^j {\Gamma}^k_{ij}(x) \frac{\partial F}{\partial y^k}=0\ \  (i=1, \ldots,n),
\end{equation}
where ${\Gamma}^k_{ij}(x)$'s are the connection parameters and the vector fields of type $\displaystyle{\partial/\partial x^i-y^j {\Gamma}^k_{ij}(x) \partial/\partial y^k}$ span the associated horizontal distribution belonging to $\nabla$. Equation (\ref{cond1}) is called the \emph{compatibility equation}. 

The concept of generalized Berwald manifolds goes back to V. Wagner \cite{Wag1}. For a summary of the recent trends and some general results in the theory of generalized Berwald manifolds see \cite{V11}. To express the compatible linear connection in terms of the canonical data of the Finsler manifold is the problem of the intrinsic characterization we are going to solve in some special cases of locally symmetric $m$-th root metrics. Especially, we are interested in generalized Berwald surfaces with locally symmetric fourth root metrics. The main result (Theorem 1) is their intrinsic characterization in terms of the basic notions of linear algebra. We present a one-parameter family of examples as well. The last section contains some computations in 3D. They are supported by the MAPLE mathematics softwer (LinearAlgebra).

\section{Locally symmetric $m$-rooth metrics}

\begin{Def}
Let $m=2l$ be a positive natural number, $l=1, 2, \ldots. $ A Finslerian metric $F$ is called an $m$-th root metric if its $m$-th power $F^m$ is of class $C^{m}$ on the tangent manifold $TM$.
\end{Def}

Using that $F$ is positively homogeneous of degree one, its $m$-th power is homogeneous of degree $m$. Since it is of class $C^m$ on the tangent manifold $TM$ (including the zero section), its local form must be a polynomial of degree $m$ in the variables $y^1$, $\ldots$, $y^n$ as follows:
\begin{equation}
\label{polgen}
F^m(x,y)=\sum_{i_1+\ldots+i_n=m} a_{i_1 \ldots i_n}(x)(y^1)^{i_1}\cdot \ldots (y^n)^{i_n}.
\end{equation}
Finslerian metrics of the form (\ref{polgen}) has been introduced by Shimada \cite{Shim}. They are generalizations of the so-called Berwald-Mo\'{o}r metrics. The geometry of the $m$-th root metrics and some special cases have been investigated by several authors such as  M. Matsumoto, K. Okubo, V. Balan, N. Brinzei, L. Tam\'{a}ssy, A. Tayebi and B. Najafi etc. in \cite{BB}, \cite{VB}, \cite{Brin}, \cite{MO}, \cite{TN} and \cite{Tam}.

\begin{Exm}
Riemannian metrics are $2$nd root metrics, i.e. $m=2$.
\end{Exm}

\begin{Def}
$F$ is locally symmetric if each point has a coordinate neighbourhood such that $F^m$ is a symmetric polynomial of degree $m$ in the variables $y^1$, $\ldots$, $y^n$ of the induced coordinate system on the tangent manifold.
\end{Def}

Suppose that formula (\ref{polgen}) is a symmetric expression of $F(x,y)$ in the variables $y^1$, $\ldots$, $y^n$. Using the fundamental theorem of symmetric polynomials, we can write that
\begin{equation}
\label{charpol}
F^m(x,y)=P(s^1, \ldots, s^n),
\end{equation}
where 
$$s^1=y^1+\ldots y^n, \ s^2=y^1y^2+\ldots+y^{n-1}y^n, \ldots, s^n=y^1\cdot \ldots \cdot y^n$$
are the so-called elementary symmetric polynomials. The polynomial $P$ with coefficients depending on the position is called the \emph{local characteristic polinomial} of the locally symmetric $m$-th root metric. Using the homogenity properties, the reduction of the number of the coefficients depending on the position is
\begin{equation}
\label{polspec}
F^m(x,y)=\sum_{j_1+2j_2+\ldots+nj_n=m} c_{j_1 \ldots j_n}(x)(s^1)^{j_1}\cdot \ldots (s^n)^{j_n}.
\end{equation}

The following tables show the possible values of the powers $j_1$, $\ldots$, $j_n$ in case of $n=2, 3, 4, 5$ and $m=4$ (fourth root metrics). The corresponding local characteristic polynomials are of the form
\begin{equation*}
\begin{aligned}
&P(s^1, s^2)=c_{40}(x)(s^1)^4+c_{21}(x)(s^1)^2s^2+c_{02}(x)(s^2)^2,\\
&P(s^1, s^2,s^3)=c_{400}(x)(s^1)^4+c_{210}(x)(s^1)^2s^2+c_{020}(x)(s^2)^2+c_{101}(x)s^1s^3,\\
&P(s^1, s^2,s^3,s^4)=c_{4000}(x)(s^1)^4+c_{2100}(x)(s^1)^2s^2+c_{0200}(x)(s^2)^2+c_{1010}(x)s^1s^3+c_{0001}(x)s^4.
\end{aligned}
\end{equation*}

\vspace{0.2cm}
\begin{center}
\begin{tabular}{|l|l|}
\hline
 \multicolumn{2}{|c|}{$n=2$}\\
\multicolumn{2}{|c|}{$j_1+2j_2=4$}\\
\hline
$j_1=4$&$j_2=0$\\
\hline
$j_1=2$&$j_2=1$\\
\hline
$j_1=0$&$j_2=2$\\
\hline
--&--\\
\hline
--&--\\
\hline
\end{tabular} \ \ \begin{tabular}{|l|l|l|}
\hline
 \multicolumn{3}{|c|}{$n=3$}\\
\multicolumn{3}{|c|}{$j_1+2j_2+3j_3=4$}\\
\hline
$j_1=4$&$j_2=0$&$j_3=0$\\
\hline
$j_1=2$&$j_2=1$&$j_3=0$\\
\hline
$j_1=0$&$j_2=2$&$j_3=0$\\
\hline
$j_1=1$&$j_2=0$&$j_3=1$\\
\hline
--&--&--\\
\hline
\end{tabular} \ \ \begin{tabular}{|l|l|l|l|}
\hline
 \multicolumn{4}{|c|}{$n=4$}\\
\multicolumn{4}{|c|}{$j_1+2j_2+3j_3+4j_4=4$}\\
\hline
$j_1=4$&$j_2=0$&$j_3=0$&$j_4=0$\\
\hline
$j_1=2$&$j_2=1$&$j_3=0$&$j_4=0$\\
\hline
$j_1=0$&$j_2=2$&$j_3=0$&$j_4=0$\\
\hline
$j_1=1$&$j_2=0$&$j_3=1$&$j_4=0$\\
\hline
$j_1=0$&$j_2=0$&$j_3=0$&$j_4=1$\\
\hline
\end{tabular}
\end{center}
\vspace{0.2cm}
The case $n=4$ also shows the special form of the local characteristic polynomial of locally symmetric fourth root metrics for $n\geq 5$ up to the formal zero powers for the terms $s_5$, $\ldots$, $s_n$. If $n=5$ then we have that
$$P(s^1, s^2,s^3,s^4,s^5)=c_{40000}(x)(s^1)^4+c_{21000}(x)(s^1)^2s^2+c_{02000}(x)(s^2)^2+c_{10100}(x)s^1s^3+c_{00010}(x)s^4.$$

\begin{Cor} A locally symmetric fourth root metric is locally determined by at most five components of its local characteristic polynomial.
\end{Cor}

\subsection{Regularity properties} Let us introduce the following notations
$$A:=F^m(x,y)=\sum_{i_1+\ldots+i_n=m} a_{i_1 \ldots i_n}(x)(y^1)^{i_1}\cdot \ldots (y^n)^{i_n}, \ \ A_i:=\frac{\partial A}{\partial y^i} \ \ \textrm{and} \ \ A_{ij}:=\frac{\partial^2 A}{\partial y^i \partial y^j}.$$
\begin{Lem} The function $F=\sqrt[m]{A}$ is a Finsler metric if and only if $A_{ij}(v)$ is a positive definite matrix for any nonzero $v \in \pi^{-1}(U)$.
\end{Lem}
\begin{Pf} 
Suppose that $A_{ij}(v)$ is a positive definite matrix for any nonzero $v \in \pi^{-1}(U)$. Then, by the homogenity properties,
$$0<y^iy^jA_{ij}=m(m-1)A,$$
i.e. we can introduce a positive valued one-homogeneous function $F=\sqrt[m]{A}$ on the complement of the zero section in $\pi^{-1}(U)$. Since $\displaystyle{A=(F^2)^l}$,
$$A_{ij}=2^llE^{l-1} \left(g_{ij}+\frac{l-1}{E}\frac{\partial E}{\partial y^i}\frac{\partial E}{\partial y^i}\right), \ \textrm{where}\ E=\frac{1}{2}F^2.$$
Taking linearly independent vertical vector fields $V_1, \ldots, V_{n-1}$ such that $V_i(v)E=0$ at a given $v\in \pi^{-1}(U)$, it follows that $V_1, \ldots, V_{n-1}, C$ is a basis of the vertical subspace $V_vTM$. We have 
$$V_k^i(v)V_k^j(v)g_{ij}(v)=\frac{V_k^i(v)V_k^j(v) A_{ij}(v)}{2^llE^{l-1}(v)}>0\ \ \textrm{and}\ \ C^iC^jg_{ij}=y^iy^jg_{ij}=2E>0$$
in the sense of Euler's theorem on homogeneous functions. This means that $g_{ij}(v)$ is a positive definite matrix for any non-zero $v$. The converse of the statement is clear because of
$$y^iy^jA_{ij}(v)=2^llE^{l-1}(v) \left(y^iy^jg_{ij}(v)+\frac{l-1}{E(v)}\left(y^i\frac{\partial E}{\partial y^i}(v)\right)^2\right),$$
i.e. if $g_{ij}(v)$ is a positive definite matrix, then $A_{ij}(v)$ is also positive definite for any non-zero $v\in TM$.
\end{Pf}

\section{Finsler surfaces with locally symmetric $4$-rooth metrics}

Let $M$ be a two-dimensional Finsler manifold (Finsler surface) with a locally symmetric fourth root metric $F=\sqrt[4]{A}$. Its local characteristic polynomial must be of the form
\begin{equation}
P(s^1, s^2)=A(x,y)=a(x) (y^1+y^2)^4 + b(x) (y^1+y^2)^2 y^1 y^2 + c(x) (y^1y^2)^2, \label{poly1}
\end{equation}
where $a(x)=c_{40}(x)$, $b(x)=c_{21}(x)$ and $c(x)=c_{02}(x)$. Differentiating (\ref{poly1})
\begin{eqnarray*}
&&A_1=\frac{\partial A}{\partial y^1}= 4 a(x) (y^1+y^2)^3 + 2 b(x) (y^1+y^2) y^1 y^2 + b(x) (y^1+y^2)^2y^2 + 2 c(x) y^1 (y^2)^2,\\
&&A_2=\frac{\partial A}{\partial y^2}= 4 a(x) (y^1+y^2)^3 + 2 b(x) (y^1+y^2) y^1 y^2 + b(x) (y^1+y^2)^2y^1 + 2 c(x) y^2 (y^1)^2.
\end{eqnarray*}
By some further computations
\begin{eqnarray}
&&A_{11}= 12 a(x) (y^1)^2 + ( 24 a(x)+ 6 b(x)) y^1 y^2 + ( 12 a(x) + 4b(x) + 2c(x) ) (y^2)^2,\label{A11}\\
&&A_{12}=A_{21}= ( 12 a(x)+ 3 b(x))(y^1)^2 + (24 a(x) + 8 b(x) + 4 c(x))  y^1 y^2 + (12 a(x) + 3b(x))(y^2)^2,\label{A12}\\
&&A_{22}= ( 12 a(x) + 4b(x) + 2c(x) )(y^1)^2 + ( 24 a(x)+ 6 b(x)) y^1 y^2 + 12 a(x) (y^2)^2 .\label{A22}
\end{eqnarray}
Introducing the functions
\begin{equation} l(x):= a(x), \ m(x):= 4 a(x) + b(x), \ n(x):= 6 a(x) + 2 b(x) + c(x),
\end{equation}
we have that 
\begin{equation}
\label{metric}
A(x,y)=l(x)(y^1)^{4}+ m(x) (y^1)^{3}y^2+n(x)(y^1)^2 (y^2)^{2}+ m(x) y^1 (y^2)^{3}+l(y^2)^{4},
\end{equation}
\begin{eqnarray}
&&\frac{\partial A}{\partial y^1}= 4 l(x) (y^1)^3 + 3 m (x) (y^1)^2 y^2+ 2n (x)y^1 (y^2)^2 + m (x) (y^2)^3, \label{A1}\\
&&\frac{\partial A}{\partial y^2}= m (x)(y^1)^3 + 2n (x) (y^1)^2 y^2+ 3m (x)y^1 (y^2)^2 + 4l (x) (y^2)^3.\label{A2}
\end{eqnarray}
Since
\[
\left[
\begin{matrix}
l \\ m \\ n
\end{matrix}
\right]
=
\left[
\begin{matrix}
1 & 0 & 0 \\ 4 & 1 & 0 \\ 6 & 2 & 1
\end{matrix}
\right]
\left[
\begin{matrix}
a \\ b \\ c
\end{matrix}
\right].
\]
is a regular linear transformation, the coefficients $a(x), b(x)$, $c(x)$ are uniquely determined by $l(x)$, $m(x)$, $n(x)$ and vice versa. Using (\ref{A11}), (\ref{A12}) and (\ref{A22}) 
$$A_{11}= \left[
\begin{matrix}
y^1 & y^2
\end{matrix}
\right]
\left[
\begin{matrix}
12l & 3m\\ 3m & 2n
\end{matrix}
\right]
\left[
\begin{matrix}
y^1 \\ y^2
\end{matrix}
\right], \ A_{12}= A_{21}=\left[
\begin{matrix}
y^1 & y^2
\end{matrix}
\right]
\left[
\begin{matrix}
3m & 2n\\ 2n & 3m
\end{matrix}
\right]
\left[
\begin{matrix}
y^1 \\ y^2
\end{matrix}
\right], $$
$$A_{22}= \left[
\begin{matrix}
y^1 & y^2
\end{matrix}
\right]
\left[
\begin{matrix}
2n & 3m\\ 3m & 12l
\end{matrix}
\right]
\left[
\begin{matrix}
y^1 \\ y^2
\end{matrix}
\right].
$$
Therefore $A_{ij}=\left[
\begin{matrix}
A_{11} & A_{12}\\ A_{21} & A_{22}
\end{matrix}
\right] $ is positive definite if and only if 
$$ \left[ \begin {array}{cc} 12\,l&3\,m\\ \noalign{\medskip}3\,m&2\,n
\end {array} \right] \ \ \textrm{and}\ \  \left[ \begin {array}{cc} 12\,l&3\,m\\ \noalign{\medskip}3\,m&2\,n
\end {array} \right] \left[ \begin {array}{cc} 2\,n&3\,m\\ \noalign{\medskip}3\,m&12\,l
\end {array} \right] - \left[ \begin {array}{cc} 3\,m&2\,n\\ \noalign{\medskip}2\,n&3\,m
\end {array} \right]^2$$
are positive definite. Using some direct computations 
$$\left[ \begin {array}{cc} 12\,l&3\,m\\ \noalign{\medskip}3\,m&2\,n
\end {array} \right] \left[ \begin {array}{cc} 2\,n&3\,m\\ \noalign{\medskip}3\,m&12\,l
\end {array} \right] - \left[ \begin {array}{cc} 3\,m&2\,n\\ \noalign{\medskip}2\,n&3\,m
\end {array} \right]^2= \left[ \begin {array}{cc} 24\,nl-4\,{n}^{2}&72\,ml-12\,nm
\\ \noalign{\medskip}0&24\,nl-4\,{n}^{2}\end {array} \right] 
$$
and, consequently, 
\begin{eqnarray}
12 l >0, \ \  24 l n - 9 m^2>0\ \ \textrm{and} \ \ 24 n l - 4 n^2 >0.\label{posi1}
\end{eqnarray}
Especially (\ref{posi1}) is equivalent to 
\begin{equation}
\label{posi2} 
6l>n>0  \ \ \textrm{and}\ \ \frac{8}{3} n l > m^2.
\end{equation}
\subsection{Generalized Berwald surfaces with locally symmetric fourth root metrics}  Let $\nabla$ be a linear connection on the base manifold $M$ equipped with a locally symmetric fourth root metric $F=\sqrt[4]{A}$ and suppose that the parallel transports preserve the Finslerian length of tangent vectors. The compatibility condition (\ref{cond1}) can be written into the form
\begin{equation}
\frac{\partial A}{\partial x^i} - y^j \Gamma^k_{ij}(x) \frac{\partial A}{\partial y^k}=0 \ \ (i=1, 2).\label{gene}
\end{equation}
Substituting (\ref{metric}), (\ref{A1}) and (\ref{A2}) into (\ref{gene}), we get the following system of linear equations 
\begin{eqnarray}
\left[
\begin{matrix}
4l & 0 & m & 0 \\ 3m & 4l & 2n & m \\ 2n & 3m & 3m & 2n \\ m & 2n & 4l & 3m \\ 0 & m & 0 & 4l
\end{matrix}
\right]
\left[
\begin{matrix}
\Gamma^1_{i1}\\ \Gamma^1_{i2} \\ \Gamma^2_{i1} \\ \Gamma^2_{i2}\\
\end{matrix}
\right]
=
\left[
\begin{matrix}
\partial l / \partial x^i \\ \partial m / \partial x^i \\ \partial n / \partial x^i \\ \partial m / \partial x^i \\ \partial l / \partial x^i
\end{matrix}
\right]\label{system}
\end{eqnarray}
because of
$$y^j \Gamma^k_{ij} \frac{\partial A}{\partial y^k}=\left(y^1 \Gamma^1_{i1}+y^2 \Gamma^1_{i2}\right) \frac{\partial A}{\partial y^1}+\left(y^1 \Gamma^2_{i1}(x)+y^2 \Gamma^2_{i2}\right) \frac{\partial A}{\partial y^2}\stackrel{(\ref{A1}), \  (\ref{A2})}{=}$$
$$$$
$$\left(4l\Gamma^1_{i1}+m\Gamma^2_{i1}\right)(y^1)^4+\left(m\Gamma^2_{i2}+2n\Gamma^2_{i1}+4l\Gamma^1_{i2}+3m\Gamma^1_{i1}\right)(y^1)^3 y^2+$$
$$$$
$$\left(2n\Gamma^2_{i2}+3m\Gamma^2_{i1}+3m\Gamma^1_{i2}+2n\Gamma^1_{i1}\right)(y^1)^2(y^2)^2+$$
$$$$
$$\left(3m\Gamma^2_{i2}+4l\Gamma^2_{i1}+2n\Gamma^1_{i2}+m\Gamma^1_{i1}\right)y^1(y^2)^3+\left(m\Gamma^1_{i2}+4l\Gamma^2_{i2}\right)(y^2)^4.$$
\begin{Lem}
\label{lemmakey}
If $F=\sqrt[4]{A}$ is a non-Riemannain connected generalized Berwald surface with a locally symmetric fourth root metric, then 
\begin{eqnarray}
\textrm{rank} \left[
\begin{matrix}
4l & 0 & m & 0 \\ 3m & 4l & 2n & m \\ 2n & 3m & 3m & 2n \\ m & 2n & 4l & 3m \\ 0 & m & 0 & 4l
\end{matrix}
\right]=4. \label{matrix}
\end{eqnarray}
\end{Lem}

\begin{Pf}
Suppose, in contrary, that the rank is less than 4. In case of $m(x)=0$ 
\begin{eqnarray*}
\textrm{rank} \left[
\begin{matrix}
4l(x) & 0 & 0 & 0 \\ 0 & 4l(x) & 2n(x) & 0 \\ 2n(x) & 0 & 0 & 2n(x) \\ 0 & 2n(x) & 4l(x) & 0 \\ 0 & 0 & 0 & 4l(x)
\end{matrix}
\right]<4.
\end{eqnarray*}
Therefore $n(x)= 2l(x)=2a(x)$. Since $m(x)=0$, it follows that $b(x)=-4a(x)$, $c(x)=4a(x)$ and, consequently, $b^2(x)-4a(x)c(x)=0$. This means that 
\begin{eqnarray*}
A(x,y)=a(x) (y^1)^4+ 2 a(x) (y^1)^2 (y^2)^2 + a(x) (y^2)^4= a(x) \left( (y^1)^2 + (y^2)^2 \right)^2
\end{eqnarray*}
is a complete square of the quadratic form $E(x,y)$ with respect to the variables $y^1$ and $y^2$. Using that we have a compatible linear connection, the indicatrices are quadrics at each point as the (linear) parallel translates of the indicatrix at the single point $x$. It is a contradiction because the surface is non-Riemannian. In the second case we suppose that $m(x) \neq 0$ and consider the submatrix 
\begin{eqnarray*}
\left[
\begin{matrix}
4l(x) & 0 & m(x) & 0 \\ 2n(x) & 3m(x) & 3m(x) & 2n(x) \\ m(x) & 2n(x) & 4l(x) & 3m(x) \\ 0 & m(x) & 0 & 4l(x)
\end{matrix}
\right].
\end{eqnarray*}
If the rank is less than 4, then its determinant must be zero:
$$4m(x) (6l(x) - n(x)) (8 l^2(x) -4 l(x) n(x) + m^2(x))=0.$$
According to the regularity properties (\ref{posi1}), $6l(x)-n(x)>0$ and $m(x)\neq 0$ imply that 
$$8 l^2(x) -4 l(x) n(x) + m^2(x)=0,$$
$$8a^2(x)-4a(x)\left(6 a(x) + 2 b(x) + c(x)\right)+\left(4 a(x) + b(x)\right)^2=0 \ \ \Rightarrow \ \ b^2(x)-4a(x)c(x)=0$$
and the proof can be finished as above. 
\end{Pf}

\begin{Thm}
\label{maintheorem2D}
Let $M$ be a connected non-Riemannian Finsler surface with a locally symmetric fourth root metric $F=\sqrt[4]{A}$. It is a generalized Berwald surface if and only if the coefficient matrix 
$$B:=\left[
\begin{matrix}
4l & 0 & m & 0 \\ 3m & 4l & 2n & m \\ 2n & 3m & 3m & 2n \\ m & 2n & 4l & 3m \\ 0 & m & 0 & 4l
\end{matrix}
\right]$$
is of constant rank $4$ and  
\begin{eqnarray}
\label{main2D}
\det \left[
\begin{matrix}
4l & 0 & m & 0 & \partial l/ \partial x^i \\ 3m & 4l & 2n & m & \partial m/ \partial x^i \\ 2n & 3m & 3m & 2n & \partial n/ \partial x^i \\ m & 2n & 4l & 3m & \partial m/ \partial x^i \\ 0 & m & 0 & 4l & \partial l/ \partial x^i
\end{matrix}
\right] =0 \ \ \ (i=1,2).
\end{eqnarray}
The compatible linear connection is uniquely determined by the formula
\begin{equation}
\label{explicite}
\left[
\begin{matrix}
\Gamma^1_{i1}\\ \Gamma^1_{i2} \\ \Gamma^2_{i1} \\ \Gamma^2_{i2}\\
\end{matrix}
\right]
=\left(B^{T} B\right)^{-1} B^T 
\left[
\begin{matrix}
\partial l / \partial x^i \\ \partial m / \partial x^i \\ \partial n / \partial x^i \\ \partial m / \partial x^i \\ \partial l / \partial x^i
\end{matrix}
\right] \ \ (i=1,2).
\end{equation}
\end{Thm}
\begin{Pf}
According to Lemma \ref{lemmakey}, if a generalized Berwald metric is non-Riemannian, then the coefficient matrix in (\ref{system}) is of constant rank $4$ and the existence of the unique solution is equivalent to the vanishing of the determinant of the extended matrix. To conclude the explicite expression of the coefficients of the compatible linear connection note that $B^TB$ is the Gram matrix of the linearly independent column vectors, i.e. the Gram determinant different from zero and $B^TB$ is invertible. Conversely, if the rank is maximal and (\ref{main2D}) holds, then the system of linear equations (\ref{system}) has a unique solution given by formula (\ref{explicite}).
\end{Pf}

\subsection{Examples} In what follows we give explicite examples for non-constant functions satisfying (\ref{main2D}). To simplify the formulation of the problem we prove that the class of generalized Berwald metrics is closed under the conformal deformation. 

\begin{Def}
The Finsler metrics $F_1$ and $F_2$ are conformally related if $F_2(x,y)=e^{\alpha(x)}F_1(x,y).$
\end{Def}

\begin{Thm}
\label{confinv}
The class of generalized Berwald metrics is closed under the conformal deformation.
\end{Thm}

\begin{Pf}
If the compatibility equation (\ref{cond1}) holds for the Finsler metric $F_1$ then we have that
$$
\frac{\partial F_2}{\partial x^i}-y^j {\Gamma}^k_{ij}(x) \frac{\partial F_2}{\partial y^k}=e^{\alpha(x)}\left(\frac{\partial F_1}{\partial x^i}-y^j {\Gamma}^k_{ij}(x) \frac{\partial F_1}{\partial y^k}\right)+F_2\frac{\partial \alpha}{\partial x^i}=F_2\frac{\partial \alpha}{\partial x^i}.
$$
This means that the compatibility equation
\begin{equation}
\label{cond2}
\frac{\partial F_2}{\partial x^i}-y^j\left({\Gamma}^k_{ij}(x)+\frac{\partial \alpha}{\partial x^i}(x)\delta_j^k\right) \frac{\partial F_2}{\partial y^k}=0\ \  (i=1, \ldots,n)
\end{equation}
holds for the Finsler metric $F_2$.
\end{Pf}

In the sense of the previous theorem we can suppose that $4l(x)=1$ (conformal deformation of the metric). Therefore we should find solutions $n(x)$ and  $m(x)$ of equation
\begin{eqnarray*}
\label{main2Dexm}
\det \left[
\begin{matrix}
1 & 0 & m & 0 & 0 \\ 3m & 1 & 2n & m & \partial m/ \partial x^i \\ 2n & 3m & 3m & 2n & \partial n/ \partial x^i \\ m & 2n &1 & 3m & \partial m/ \partial x^i \\ 0 & m & 0 & 1 & 0
\end{matrix}
\right] =0 \,\,\,\,\,\,\,\,\,\,\,\, (i=1,2).
\end{eqnarray*}
We have 
$$8\partial m/ \partial x^i {m}^{3}n-8\,\partial n/ \partial x^i {m}^{4}-12\partial m/ \partial x^i {m}^{3}-8\partial m/ \partial x^i m {n}^{2}+12\partial n/ \partial x^i {m}^{2}n+$$
$$16\partial m/\partial x^i mn-2\partial n/ \partial x^i{m}^{2}-4\partial n/ \partial x^i {n}^{2}-6m\partial m/ \partial x^i +\partial n/ \partial x^i=0,$$
\begin{equation}
\label{exm1}
2m \left( 2\,n-3 \right)  \left( 2\,{m}^{2}-2\,n+1 \right) \partial m/ \partial x^i+ \left( 2
\,{m}^{2}-2\,n+1 \right)  \left( -4\,{m}^{2}+2\,n+1 \right) \partial n/ \partial x^i=0.
\end{equation}
Since $l=1/4$, the case $2{m}^{2}-2n+1=0$ gives that $m^2=n-1/2$,
$$8 l^2(x) -4 l(x) n(x) + m^2(x)=0,$$
$$8a^2(x)-4a(x)\left(6 a(x) + 2 b(x) + c(x)\right)+\left(4 a(x) + b(x)\right)^2=0 \ \ \Rightarrow \ \ b^2(x)-4a(x)c(x)=0,$$
i.e. the surface is Riemannian\footnote{Recall that the quadratic indicatrix at a single point implies that the indicatrices are quadratic at all points of the manifold because of the compatible (linear) parallel transports.}. Therefore we can suppose that $2{m}^{2}(x)-2n(x)+1\neq 0$ to present non-Riemannian generalized Berwald surfaces. We have 
$$2m \left( 2n-3 \right) \partial m/ \partial x^i+ \left( -4{m}^{2}+2n+1 \right) \partial n/ \partial x^i=0,$$
$$2m \partial m/ \partial x^i=-{\frac { -4{m}^{2}+2n+1}{2n-3}}\partial n/ \partial x^i.$$
Recall that the regularity conditions in (\ref{posi2}) reduce to
\begin{equation}
\label{posi3}
\frac{3}{2}>n>0 \ \ \textrm{and}\ \ \frac{2}{3}> \frac{m^2}{n}
\end{equation}
because of $l=1/4$ and we can divide by the term $2n-3$. To present non-constant\footnote{The simplest examples are the constant functions satisfying (\ref{matrix}). In this case (\ref{main2D}) is automathic.} solutions $m(x)$ and $n(x)$ suppose that $n$ is regular at the point $x$ and, by taking $m^2=y(n)$, we have that 
$$y'(n) =-{\frac {-4y \left( n \right) +2n+1}{2n-3}}.$$
The general form of the solution is
$$y_{\varepsilon} \left( n \right) =n-1/2+ \varepsilon \left( 2\,n-3 \right) ^{2}.$$
The integration constant $\varepsilon=1/18$ provides that it is a positive-valued function on the positive half line of the reals. On the other hand 
$$\frac{m_{1/18}^2}{n}=\frac{y_{1/18}(n)}{n}=\frac{1}{3}+\frac{2}{9}n < \frac{2}{3}$$ 
on the interval $(0, 3/2)$. Therefore the regularity conditions in (\ref{posi3}) are satisfied  and  the triplet
$$l(x)=r(x)/4, \ \ m_{1/18}(x):=r(x) \sqrt{n-1/2+ \frac{1}{18} \left( 2\,n-3 \right) ^{2}},\ \ r(x) n(x)$$
determines a non-Riemannian generalized Berwald surface for any positive-valued function $r(x)$; see Theorem \ref{maintheorem2D} and Theorem \ref{confinv}. It is a one-parameter family of examples depending on the integration constant $\varepsilon$ satisfying the regularity conditions in (\ref{posi3}). For example for any $1/18 < \varepsilon < 1/6$ it follows that
$$0 < n_0:=\frac{18\varepsilon-1}{8\varepsilon}< \frac{3}{2}$$
and 
$$\frac{m_{\varepsilon}^2(n_0)}{n_0}=\frac{y_{\varepsilon}(n_0)}{n_0}=\varepsilon+\frac{1}{2}< \frac{1}{6}+\frac{1}{2}= \frac{2}{3}.$$
Using a continuity argument, if $n$ is sufficiently close to $n_0$, then the  the regularity conditions in (\ref{posi3}) are satisfied. 

\section{Computations in 3D}

Suppose that the base manifold is of dimension $3$, i.e. the locally symmetric fourth root metric must be of the form 
$$P(s^1, s^2, s^3)=A(x,y)=a(x) (y^1+y^2+y^3)^4 + b(x) (y^1+y^2+y^3)^2 ( y^1 y^2 + y^1 y^3 + y^2 y^3 ) +$$
$$ c(x) ( y^1 y^2 + y^1 y^3 + y^2 y^3 )^2 +d(x) (y^1+y^2+y^3)y^1 y^2 y^3,$$
where $a(x)=c_{400}(x)$, $b(x)=c_{210}(x)$, $c(x)=c_{020}(x)$ and $d(x)=c_{101}(x)$. Introducing the functions
$$l(x):=a(x), \ m(x):= 4a(x)+b(x), \ n(x):= 6a(x) + 2b(x) +c(x),$$
$$q(x) := 12 a(x)+ 5b(x) + 2c(x)+d(x),$$
we have that 
$$
A=l \left((y^1)^4+(y^2)^4+(y^3)^4 \right)+m\left((y^1)^3 y^2 + (y^1)^3 y^3 + (y^2)^3 y^3 + (y^3)^3 y^2 + (y^3)^3 y^1 + (y^2)^3 y^1 \right)+$$
$$n\left((y^1)^2( y^2)^2 + (y^1)^2( y^3)^2 + (y^2)^2( y^3)^2\right)+q\left((y^1)^2y^2y^3+y^1( y^2)^2 y^3+ y^1 y^2( y^3)^2 \right).$$
By some further computations
\begin{eqnarray*}
A_1:= \frac{\partial A}{\partial y^1}&=& 4l(y^1)^3+m \left(3(y^1)^2 y^2+3(y^1)^2 y^3+(y^2)^3+(y^3)^3\right)+n\left(2y^1 (y^2)^2+2y^1(y^3)^2\right)+\\
&&q\left(2y^1y^2y^3+(y^2)^2y^3+y^2(y^ 3)^2\right),\\
A_2:= \frac{\partial A}{\partial y^2}&=& 4l(y^2)^3+m\left((y^1)^3+3y^1(y^2)^2+3(y^2)^2y^3+(y^3)^3\right)+n\left(2(y^1)^2y^2+2y^2(y^3)^2\right)+\\
&&q\left((y^1)^2y^3+2y^1y^2y^3+y^1(y^3)^2\right),\\
A_3:= \frac{\partial A}{\partial y^3}&=& 4l(y^3)^3+m\left((y^1)^3+3y^1(y^3)^2+(y^2)^3+3y^2(y^3)^2\right)+n\left(2(y^1)^2 y^3+2(y^2)^2y^3\right)+\\
&&q\left((y^1)^2y^2+y^1(y^2)^2+2y^1y^2y^3\right),\\
A_{11}&=&  12 l (y^1)^2 + m \left( 6y^1y^2+6y^1y^3\right)+ n \left(2(y^2)^2 + 2(y^3)^2\right) + 2q y^2y^3,\\
A_{12}&=& A_{21}=m \left( 3(y^1)^2 + 3(y^2)^2\right)+ 4n y^1y^2 + q \left(2y^1y^3+ 2y^2y^3 + (y^3)^2\right),\\
A_{13}&=&A_{31}=  m \left( 3(y^1)^2+3(y^3)^2 \right)+4n y^1y^3+ q \left( 2y^1y^2 + (y^2)^2 + 2y^2y^3\right) ,\\
A_{22}&=& 12 l (y^2)^2 + m \left(6 y^1y^2+ 6 y^2 y^3\right)+ n \left(2 (y^1)^2 + 2 (y^3)^2\right)+ 2q  y^1y^3,\\
A_{23}&=&A_{32}= m \left( 3(y^2)^2 + 3(y^3)^2\right)+ 4n y^2y^3  + q \left((y^1)^2 + 2y^1y^2 + 2y^1y^3\right) ,\\
A_{33}&=& 12 l (y^3)^2 + m \left(6 y^1y^3+ 6 y^2y^3\right)+ n \left(2 (y^1)^2 + 2 (y^2)^2\right)+ 2q y^1y^2.
\end{eqnarray*}
The second order partial derivatives can be written into the form 
$$A_{11}=
\left[
\begin{matrix}
y^1 & y^2 & y^3 
\end{matrix}
\right] 
\left[
\begin{matrix}
12 l & 3m & 3m \\  3m & 2n & q \\  3m & q & 2n
\end{matrix}
\right]
\left[
\begin{matrix}
y^1 \\  y^2 \\ y^3
\end{matrix}
\right], \ A_{12}=A_{21}=
\left[
\begin{matrix}
y^1 & y^2 & y^3 
\end{matrix}
\right]
\left[
\begin{matrix}
3m & 2n & q \\  2n & 3m & q \\  q & q & q
\end{matrix}
\right]
\left[
\begin{matrix}
y^1 \\  y^2 \\ y^3
\end{matrix}
\right],$$
$$A_{22}=
\left[
\begin{matrix}
y^1 & y^2 & y^3 
\end{matrix}
\right]
\left[
\begin{matrix}
2n & 3m & q \\  3m & 12l & 3m \\  q & 3m & 2n
\end{matrix}
\right]
\left[
\begin{matrix}
y^1 \\  y^2 \\ y^3
\end{matrix}
\right], \ A_{13}=A_{31}=
\left[
\begin{matrix}
y^1 & y^2 & y^3 
\end{matrix}
\right]
\left[
\begin{matrix}
3m & q & 2n \\  q & q & q \\  2n & q & 3m
\end{matrix}
\right]
\left[
\begin{matrix}
y^1 \\  y^2 \\ y^3
\end{matrix}
\right],$$
$$ A_{33}=
\left[
\begin{matrix}
y^1 & y^2 & y^3 
\end{matrix}
\right]
\left[
\begin{matrix}
2n & q & 3m \\  q & 2n & 3m \\  3m & 3m & 12l
\end{matrix}
\right]
\left[
\begin{matrix}
y^1 \\  y^2 \\ y^3
\end{matrix}
\right], \ A_{23}=A_{32}= 
\left[
\begin{matrix}
y^1 & y^2 & y^3 
\end{matrix}
\right]
\left[
\begin{matrix}
q & q & q \\  q & 3m & 2n \\  q & 2n & 3m
\end{matrix}
\right]
\left[
\begin{matrix}
y^1 \\  y^2 \\ y^3
\end{matrix}
\right]. $$
The positive definiteness of $\displaystyle{A_{ij}:=
\left[
\begin{matrix}
A_{11} & A_{12} & A_{13} \\  A_{21} & A_{22} & A_{23} \\  A_{31} & A_{32} & A_{33}
\end{matrix}
\right]}$ can be expressed in terms of the formal subdeterminants
$$\left[
\begin{matrix}
12 l & 3m & 3m \\  3m & 2n & q \\  3m & q & 2n
\end{matrix}
\right], \ \left[
\begin{matrix}
12 l & 3m & 3m \\  3m & 2n & q \\  3m & q & 2n
\end{matrix}
\right]\left[
\begin{matrix}
2n & 3m & q \\  3m & 12l & 3m \\  q & 3m & 2n
\end{matrix}
\right]-\left[
\begin{matrix}
3m & 2n & q \\  2n & 3m & q \\  q & q & q
\end{matrix}
\right]^2, \ \ldots
$$
of the matrix (\ref{bigmatrix}) all of whose elements are matrices. Therefore the formal determinants are also matrices and all of them must be positive definite. 
\begin{equation}
\label{bigmatrix}
\left[
\begin{matrix}
\left[
\begin{matrix}
12 l & 3m & 3m \\  3m & 2n & q \\  3m & q & 2n
\end{matrix}
\right] & \left[
\begin{matrix}
3m & 2n & q \\  2n & 3m & q \\  q & q & q
\end{matrix}
\right] & \left[
\begin{matrix}
3m & q & 2n \\  q & q & q \\  2n & q & 3m
\end{matrix}
\right] \\  
&& \\ \left[
\begin{matrix}
3m & 2n & q \\  2n & 3m & q \\  q & q & q
\end{matrix}
\right] & \left[
\begin{matrix}
2n & 3m & q \\  3m & 12l & 3m \\  q & 3m & 2n
\end{matrix}
\right] & \left[
\begin{matrix}
q & q & q \\  q & 3m & 2n \\  q & 2n & 3m
\end{matrix}
\right] \\ 
&& \\ \left[
\begin{matrix}
3m & q & 2n \\  q & q & q \\  2n & q & 3m
\end{matrix}
\right] & \left[
\begin{matrix}
q & q & q \\  q & 3m & 2n \\  q & 2n & 3m
\end{matrix}
\right] & \left[
\begin{matrix}
2n & q & 3m \\  q & 2n & 3m \\  3m & 3m & 12l
\end{matrix}
\right]
\end{matrix}
\right].
\end{equation}
The compatibility equation (\ref{cond1}) can be written into the form of the following system of linear equations:
\begin{equation*}
\frac{\partial A}{\partial x^i} - y^j \Gamma^k_{ij} \circ \pi \frac{\partial A}{\partial y^k}=0 \ \ (i=1, 2),\label{gene3D}
\end{equation*}
\begin{eqnarray}
\left[
\begin{matrix}
4l & 0 & 0 & m & 0 & 0 & m & 0 & 0 \\ 0 & m & 0 & 0 & 4l & 0 & 0 & m & 0 \\ 0 & 0 & m & 0 & 0 & m & 0 & 0 & 4l \\ 3m & 4l & 0 & 2n & m & 0 & q & m & 0 \\ 3m & 0 & 4l & q & 0 & m & 2n & 0 & m \\ 0 & q & m & 0 & 3m & 4l & 0 & 2n & m \\ 0 & m & q & 0 & m & 2n & 0 & 4l & 3m \\ m & 0 & 2n & m & 0 & q & 4l & 0 & 3m \\ m & 2n & 0 & 4l & 3m & 0 & m & q & 0 \\ 2q & 3m & 3m & 2q & q & 2n & 2q & 2n & q \\ q & 2q & 2n & 3m & 2q & 3m & 2n & 2q & q \\ q & 2n & 2q & 2n & q & 2q & 3m & 3m & 2q \\ 0 & q & q & 0 & 2n & 3m & 0 & 3m & 2n \\ 2n & 0 & 3m & q & 0 & q & 3m & 0 & 2n \\ 2n & 3m & 0 & 3m & 2n & 0 & q & q & 0
\end{matrix}
\right]
\left[
\begin{matrix}
\Gamma^1_{i1} \\ \Gamma^1_{i2} \\ \Gamma^1_{i3} \\ \Gamma^2_{i1} \\ \Gamma^2_{i2} \\ \Gamma^2_{i3} \\ \Gamma^3_{i1} \\ \Gamma^3_{i2} \\ \Gamma^3_{i3}
\end{matrix}
\right]= \left[
\begin{matrix}
\partial l / \partial x^i \\ \partial l / \partial x^i \\ \partial l / \partial x^i \\ \partial m / \partial x^i \\ \partial m / \partial x^i \\ \partial m / \partial x^i \\ \partial m / \partial x^i \\ \partial m / \partial x^i \\ \partial m / \partial x^i \\  \partial q / \partial x^i \\ \partial q / \partial x^i \\ \partial q / \partial x^i \\  \partial n / \partial x^i \\ \partial n / \partial x^i \\ \partial n / \partial x^i \\
\end{matrix}
\right]. \label{bigsystem}
\end{eqnarray}

\subsection{An example} Suppose\footnote{It is also possible but the symbolical Maple computations give extremely long formulas in general.} that $4l(x)=1$, $3m(x)=1$ and $n(x)=q(x)$. According to (\ref{bigmatrix}), the matrix $A_{ij}$ is positive definite if and only if 
\begin{itemize}
\item[(i)] $\left[ \begin {array}{ccc} 3&1&1\\ \noalign{\medskip}1&2\,n&n
\\ \noalign{\medskip}1&n&2\,n\end {array} \right] 
$ is positive definite: $6n-1>0$ and $9n^2-2n>0$,
\item[(ii)] $ \left[ \begin {array}{ccc} 3&1&1\\ \noalign{\medskip}1&2\,n&n
\\ \noalign{\medskip}1&n&2\,n\end {array} \right] 
\left[ \begin {array}{ccc} 2\,n&1&n\\ \noalign{\medskip}1&3&1
\\ \noalign{\medskip}n&1&2\,n\end {array} \right] -\left[ \begin {array}{ccc} 1&2\,n&n\\ \noalign{\medskip}2\,n&1&n
\\ \noalign{\medskip}n&n&n\end {array} \right]^2=$
$$\left[ \begin {array}{ccc} -5\,{n}^{2}+7\,n&-{n}^{2}-4\,n+7&-3\,{n}^{
2}+4\,n+1\\ \noalign{\medskip}0&-5\,{n}^{2}+7\,n&-{n}^{2}+2\,n
\\ \noalign{\medskip}-{n}^{2}+2\,n&-3\,{n}^{2}+4\,n+1&{n}^{2}+2\,n
\end {array} \right] $$
is positive definite: $-5n^2+7n >0$ and $6n^3(9n^3-27n^2+23n-4) >0$, 
\item[(iii)] following the formal Sarrus rule, 
$$\left[ \begin {array}{ccc} 3&1&1\\ \noalign{\medskip}1&2\,n&n
\\ \noalign{\medskip}1&n&2\,n\end {array} \right] 
\left[ \begin {array}{ccc} 2\,n&1&n\\ \noalign{\medskip}1&3&1
\\ \noalign{\medskip}n&1&2\,n\end {array} \right]  \left[ \begin {array}{ccc} 2\,n&n&1\\ \noalign{\medskip}n&2\,n&1
\\ \noalign{\medskip}1&1&3\end {array} \right] +...=$$
$$ \left[ \begin {array}{ccc} -5\,{n}^{2}+7\,n+1&-3\,{n}^{3}-10\,{n}^{2}
+21\,n+1&-6\,{n}^{3}-9\,{n}^{2}+23\,n+11\\ \noalign{\medskip}6\,{n}^{3
}-8\,{n}^{2}-4\,n-1&-9\,{n}^{3}+7\,{n}^{2}+n&-3\,{n}^{3}-10\,{n}^{2}+
21\,n+1\\ \noalign{\medskip}18\,{n}^{3}-18\,{n}^{2}-14\,n-2&6\,{n}^{3}
-8\,{n}^{2}-4\,n-1&-5\,{n}^{2}+7\,n+1\end {array} \right] 
$$
is positive definite: $-5n^2+7n+1>0$, 
$$18\,{n}^{6}+81\,{n}^{5}-316\,{n}^{4}+154\,{n}^{3}+96\,{n}^{2}+26\,n+1>0,$$
$$-1026\,{n}^{9}+1260\,{n}^{8}+6015\,{n}^{7}-12069\,{n}^{6}+4524\,{n}^{5
}+2509\,{n}^{4}-1582\,{n}^{3}-42\,{n}^{2}+100\,n+11>0.
$$
\end{itemize}
Items (i) and (ii) give that $\displaystyle{4/3> n > \frac{5-\sqrt{13}}{6}}$. Item (ii) obviously implies that $-5n^2+7n+1>0$ in item (iii). Figure 1 shows
$$18\,{n}^{6}+81\,{n}^{5}-316\,{n}^{4}+154\,{n}^{3}+96\,{n}^{2}+26\,n+1,$$
$$-1026\,{n}^{9}+1260\,{n}^{8}+6015\,{n}^{7}-12069\,{n}^{6}+4524\,{n}^{5}+2509\,{n}^{4}-1582\,{n}^{3}-42\,{n}^{2}+100\,n+11$$
as the functions of the variable $n$ on the intervall $\displaystyle{\left(\frac{5-\sqrt{13}}{6},4/3\right)}$:
\begin{verbatim}
with(plots);
F := plot(18*n^6+..., n = 5/6-(1/6)*sqrt(13) .. 4/3, linestyle = dashdot):
G := plot(-1026*n^9+..., n = 5/6-(1/6)*sqrt(13) .. 4/3):
display({F, G});
\end{verbatim}
\begin{figure}
\includegraphics[scale=0.45]{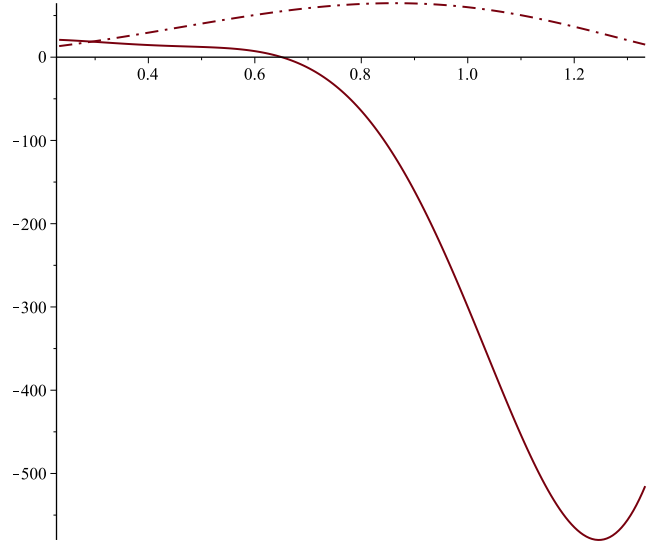}
\caption{Regularity conditions.}
\end{figure}
Therefore the regularity condition for $n$ is to be in $\displaystyle{\left(\frac{5-\sqrt{13}}{6}, r_1\right)}$, where $r_1$ is the first positive root of the polynomial $-1026\,{n}^{9}+1260\,{n}^{8}+6015\,{n}^{7}-12069\,{n}^{6}+\ldots$ For some numerical estimations note that $r_1 > 0.6458$ and $\displaystyle{\frac{5-\sqrt{13}}{6}< 0.2325}$. 
To finish the discussion of the problem consider the coefficient matrix $B$ of the system (\ref{bigsystem}). It is of rank $9$ if and only if $\det B^T B \neq 0$. Since 
$$\det B^TB = \frac{1}{387420489} (46656 n^4-5184 n^3-5589 n^2+2124 n+1492 )$$
$$(3n-2)^2(124659 n^6-51030 n^5-77517 n^4-51030 n^3+158931 n^2-92820 n+17731)^2$$
it follows that the regularity condition $0.2325 < n < 0.6458$ provides the matrix $B^TB$ to be invertible as Figure 2 shows:
\begin{verbatim}
with(plots); 
K := plot(46656*n^4..., n = 5/6-(1/6)*sqrt(13) .. 1.2, linestyle = dashdot);
L := plot(124659*n^6..., n = 5/6-(1/6)*sqrt(13) .. 1.2);
display({K, L});
\end{verbatim}
\begin{figure}
\includegraphics[scale=0.45]{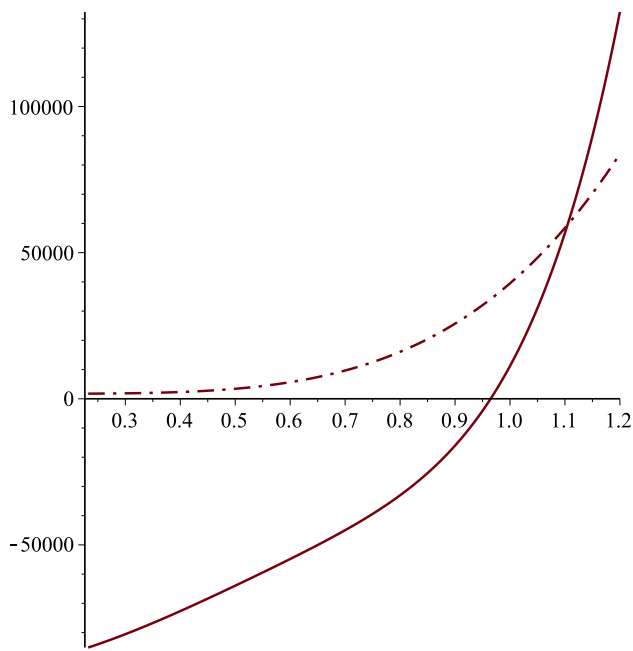}
\caption{The maximality of the rank.}
\end{figure}
Therefore the necessary and sufficient condition for the metric $F=\sqrt[4]{A}$ ($4l=1, 3m=1, n=q$, $0.2325 < n < 0.6458$ ) to be a non-Riemannian generalized Berwald metric is the vanishing of the determinant of the extended matrix of the system (\ref{bigsystem}). On the other hand, the compatible linear connection is uniquely determined by the three-dimessional analogue of formula (\ref{explicite}). Having no any other information about the linear connection, the unicity is not an automathic consequence in 3D.

\end{document}